\documentclass{amsart}

\usepackage{amssymb}
\usepackage{amsmath}

\def\wid{\operatorname{wid}}

\begin{document}

\theoremstyle{plain}
\newtheorem{theorem}{Theorem}

\markboth{Yu.~V.~Sosnovsky}
{On the Width of Verbal Subgroups}

\title[On the Width of Verbal Subgroups]
{ON THE WIDTH OF VERBAL SUBGROUPS\\
OF THE GROUPS OF TRIANGULAR MATRICES\\
OVER A~FIELD OF ARBITRARY CHARACTERISTIC}

\author{Yu.~V.~Sosnovsky}

\address{Sosnovsky Yury Vasil\/$'$evich,
Novosibirsk State Pedagogical University,
28 Vilyu\u\i{}skaya street, Novosibirsk 630126 RUSSIA}

\email{Yury\_sosnovsky@mail.ru}

\maketitle

\begin{abstract}
The width
$\wid(G,W)$
of the verbal subgroup
$v(G,W)$
of a~group~%
$G$
defined by a~collection of group words~%
$W$
is the smallest number~%
$m$
in
$\mathbb N \cup \{ +\infty \}$
such that
every element of
$v(G,W)$
is %can be represented as
the product of at most~%
$m$
words in~%
$W$
evaluated on~% the group
$G$
and their inverses.

Recall that
every verbal subgroup of the group
$T_{n} (K)$
of triangular matrices over an~arbitrary field~$K$
can be defined by just one word:
an~outer commutator word or a~power word.
We prove that
for every outer commutator word~%
$w$
the equality
$\wid(T_{n} (K),w)=1$
holds on %the group
$T_{n} (K)$
and that
if
$w=x^{s}$
then
$\wid(T_{n}(K),w)=1$
except in two cases:
\begin{enumerate}
\item[(1)]
the field~%
$K$
is finite
and~%
$s$
is divisible by the characteristic~%
$p$
of~%
$K$
but not by
$\left|K\right|-1$;

\item[(2)]
the field~%
$K$
is finite
and
$s=p^{t} (\left|K\right|-1)^{u} r$
for
$r,t,u\in \mathbb N$
with
$n\ge p^{t} +3$,
while~%
$r$
not divisible by~%
$p$.
\end{enumerate}
In these cases the width equals~2.
For finitary triangular groups
the situation is similar,
but in the second case the restriction
$n\ge p^{t} +3$
is superfluous.
\end{abstract}


\section{Introduction and Main Results}

Verbal subgroups are widely used in group theory.
Given a~group~%
$G$,
we may regard the verbal subgroup
$v(G,W)$
defined by a~certain collection of group words~%
$W$
as a~measure of distance from~%the group
$G$
to the groups of the variety
described by the collection~%of words
$W$.
It is equally important that
verbal subgroups are endomorphically admissible,
and so they are widely used as we pass to quotient groups.
To describe the set of all verbal subgroups of a~known class of groups
is a~current problem.

A~complete description of the set of verbal subgroups of
the groups
$T_{n} (K)$,
$FT(K)$,
$UT_{n} (K)$,
and
$FUT(K)$
of upper triangular,
finitary upper triangular,
upper unitriangular,
and finitary upper unitriangular
matrices over an~arbitrary field
$K$
of characteristic~0
is obtained in~\cite{4}.
It turns out that
the verbal subgroups of %the groups
$UT_{n} (K)$
and
$FUT(K)$
are precisely the members of the lower central series of these groups.
In %the groups
$T_{n} (K)$
and
$FT(K)$
the list includes in addition the subgroups
$UT_{n} (K)(D_{n} (K))^{s}$
and
$FUT(K)(FD(K))^{s}$
for
$s\in \mathbb N$
respectively.
A~similar result for fields of arbitrary characteristic
is announced in~\cite{13} and proved in~\cite{5}.
The description of verbal subgroups implies that
every verbal subgroup in these groups
is defined by just one word:
a~power word or an~outer commutator word.
Recall that
an~outer commutator word
is a~group word obtained from a~finite sequence of distinct variables
$x_{1} ,x_{2} ,\dots,x_{r}$
using only group commutators.
For instance,
$[[x_{1} ,x_{2} ],[x_{3} ,x_{4} ]]$
is an~outer commutator word,
while
$[x_{1} ,x_{2} ][x_{3} ,x_{4} ]$
and
$[[x_{1} ,x_{2} ],[x_{1} ,x_{3} ]]$
are not
since the first of them uses group multiplication
along with commutators
and the second is obtained using commutators
from the sequence
$x_{1} ,x_{2} ,x_{1} ,x_{3}$.
Recall also that
the groups of finitary matrices
$FT(K)$
($FUT(K)$)
consist of infinite invertible triangular (unitriangular) matrices
whose entries are indexed by positive integers
and differ from the entries of the identity matrix
only in finitely many entries.

Once a~description of all verbal subgroups of a~given group is available,
it is natural to become interested in their widths.
The concept of width of a~verbal subgroup,
which we recalled in the abstract,
was introduced by Merzlyakov~[6]
and turned out very useful.
Some questions immediately arise
concerning the width of a~verbal subgroup:
is it finite or not,
and if it is finite
then what is it equal to?
Quite a~few interesting results are presently available
in this subject,
in particular,
concerning the width of the commutant with respect to the commutator
and the use of these results
to characterize finite index subgroups
of finitely generated profinite groups
(for instance,
see the articles [1--3, 8--11, 14] and the book~\cite{14}).
As for open questions,
we mention the well-known question of Kargapolov [7, Question~4.34]:
with respect to which words~%
$w$
is the width
$\wid(G,w)$
of the verbal subgroup
$v(G,w)$
of a~finitely generated solvable group~%
$G$
finite?
Is
$w=x^{-1} y^{-1} xy$
among them?

Among the verbal subgroups of finite width,
the verbal subgroups of width~1
usually enjoy certain additional properties.
For instance,
%the equality
$\wid(G,x^{s} )=1$
implies that
$v(G,x^{st} )=v(v(G,x^{s} ),x^{t} )$
and that
we can extract the root of degree~%
$s$
from every element of the verbal subgroup
$v(G,x^{s} )$.
The equality
$\wid(G,\; [x,y])=1$
means that
every element of the commutant of~%the group
$G$
is a~commutator,
but this is not so when
$\wid(G,[x,y])>1$.

It is shown in~\cite{5} that
$\wid(UT_{n} (K),w)=1$
when~%
$w$
is an~outer commutator word
or
$w=x^{s}$
for
$s\in \mathbb N$
%and
%$s$
%which is
not divisible by the characteristic~%
$p$
of the field~%
$K$.
However,
if
$w=x^{s} $,
where
$s=p^{t} r$
with
$r,t\in \mathbb N$
and~%
$r$
is not divisible by~%
$p$,
then
$\wid(UT_{n} (K),x^{s} )=1$
for
$n\le p^{r} +2$
and
$\wid(UT_{n} (K),x^{s} )=2$
for
$n\ge p^{r} +3$.
For the group
$FUT(K)$
the situation regarding the width of verbal subgroups
is even simpler.
Indeed,
$\wid(FUT(K),w)=1$
when~%
$w$
is an~outer commutator word
or
$w=x^{s}$
and
$s$
is not divisible by the characteristic of~% the field~%
$K$.
However,
if~%
$s$
is divisible by the characteristic of~%the field
$K$
then
$\wid(FUT(K),x^{s} )=2$.

In this article we prove the following theorems.

\begin{theorem}
Take a~field~%
$K$
of characteristic
$p\ge 0$
with at least three elements
and arbitrary integers
$n\ge 2$
and
$s\ge 1$.
The width of the verbal subgroup of the group
$T_{n} (K)$
of triangular matrices
defined by the power word
$x^{s}$
satisfies the following equalities:
\begin{enumerate}
\item[(a)]
$wid_{} (T_{n} (K),x^{s} )=1$
when~%
$s$
is not divisible by~%
$p$;
\item[(b)]
$wid_{} (T_{n} (K),x^{s} )=2$
when~%
$s$
is divisible by~%
$p$
and
%$s$
%is
not divisible by~%
$\left|K\right|-1$
in the case that~% the field~%
$K$
is finite;
\item[(c)]
if~% the field~%
$K$
is finite
and
$s=p^{t} (\left|K\right|-1)^{u} r$
for some
$r,t,u\in \mathbb N$,
where~%
$r$
is coprime to~%
$p$
and to
$\left|K\right|-1$,
then
$wid_{} (T_{n} (K),x^{s} )=1$
for
$n\le p^{t} +2$
and
$\wid(T_{n} (K),x^{s} )=2$
for
$n\ge p^{t} +3$.
\end{enumerate}
\end{theorem}

\begin{theorem}
Take a~field~%
$K$
of characteristic
$p\ge 0$
with at least three elements
and a~positive integer~%
$s$.
The width of the verbal subgroup of
the group
$FT(K)$
of finitary triangular matrices
defined by the power word
$x^{s}$
satisfies the following equalities:
\begin{enumerate}
\item[(a)]
$\wid(FT(K),x^{s} )=1$
when~%
$s$
is not divisible by~%
$p$;
\item[(b)]
$\wid(FT(K),x^{s} )=2$
when~%
$s$
is divisible by~%
$p$.
\end{enumerate}
\end{theorem}

\begin{theorem}
The equality
$\wid(T_{n} (K),w)=1$
holds in the group
$T_{n} (K)$
of triangular matrices
for every field~%
$K$,
every outer commutator word~%
$w$,
and every integer
$n\ge 2$.
\end{theorem}

\begin{theorem}
The equality
$\wid(FT(K),w)=1$
holds in the group
$FT(K)$
of finitary triangular matrices
for every field~%
$K$
and every outer commutator word~%
$w$.
\end{theorem}

The description of verbal subgroups in %the groups
$T_{n} (K)$
and
$FT(K)$
implies that
all of them can be defined by just one word~%
$w$,
which is an~outer commutator word or a~power word.
Thus,
Theorems 1--4 enable us to compute
the width of all possible verbal subgroups of %the groups
$T_{n} (K)$
and
$FT(K)$
with respect to the natural sets of generators:
the values of~% the word
$w$.
Observe that
in case~(c) of Theorem~1
the verbal subgroup
$v(T_{n} (K),x^{s} )$
can be defined also by the outer commutator word
$[[x_{1} ,x_{2} ],[x_{3} ,x_{4} ],\dots,[x_{2p^{t} -1} ,x_{2p^{t} } ]]$;
furthermore,
the width of this verbal subgroup with respect to this word
always equals~1
independently of~%
$n$.

\section{Proofs of the Theorems}

\begin{proof}[Proof of Theorem~1]
We prove claim~(a) of Theorem~1 by induction on~%
$n$.
For
$n=2$,
given elements
$\alpha$,
$\beta$,
$\gamma$,
and
$\delta$
of~%
$K$
with
$\alpha \ne 0$
and
$\gamma \ne 0$,
the equality
\[
\left(
\begin{array}{cc}
{\alpha } & {\beta }
\\
{0} & {\gamma } \end{array}
\right)^{s}
=
\left(
\begin{array}{cc}
{\alpha ^{s} } &
  {(\alpha ^{s-1} +\alpha ^{s-2} \gamma +\dots+\alpha \gamma ^{s-2}
    +\gamma ^{s-1} )\beta }
\\
{0} & {\gamma ^{s} }
\end{array}
\right)
\]
shows that
we can obtain every matrix
\[
\left(
\begin{array}{cc}
{\alpha ^{s} } & {\delta }
\\
{0} & {\gamma ^{s} }
\end{array}
\right)
\]
in
$v(T_{2} (K),x^{s} )$
by raising to power~%
$s$
some matrix in
$T_{2} (K)$
found by choosing suitable~%
$\beta$
provided that
$\lambda =\alpha ^{s-1} +\alpha ^{s-2} \gamma +\dots+\alpha \gamma ^{s-2} +\gamma ^{s-1}
$
is nonzero.
However,
if
$\lambda =0$
then
$\alpha ^{s} =\gamma ^{s}$
and
$\alpha/\gamma$ % рекомендация автора ТеХа вместо
%${\tfrac{\alpha }{\gamma }}$
is a~degree~%
$s$
root of unity.
Replacing~%
$\gamma$
with
$\alpha$,
we arrive at %the equality
$\lambda =s\alpha ^{s-1}$
and %the inequality
$\lambda \ne 0$
since~%
$s$
is not divisible by~%
$p$.

Suppose that
we can obtain every matrix~%
$a$
in
$v(T_{n} (K),x^{s} )$
by raising to power~%
$s$
some matrix~%
$c$
in
$T_{n} (K)$
and that
either the diagonal entries of~%the matrix
$c$
are equal
or their ratios are not degree~%
$s$
roots of unity.
Represent an~arbitrary matrix~%
$\bar{a}$
in
$T_{n+1} (K)$
as
$\left(
  \begin{smallmatrix}
    {a} & {b} \\ {0} & {\gamma }
  \end{smallmatrix}
\right)$
for some matrix~%
$a$
in
$T_{n} (K)$,
a~column vector~%
$b$
of height~%
$n$,
and an~element~%
$\gamma$
of~%the field
$K$.
If~%the matrix
$\bar{a}$
lies in
$v(T_{n+1} (K),x^{s} )$
then~%the matrix
$a$
lies in
$v(T_{n} (K),x^{s} )$
and
$\gamma =\delta ^{s}$
for some~%
$\delta$
in~%
$K$.
By the inductive assumption,
there exists a~matrix~%
$c$
in
$T_{n} (K)$
such that
$a=c^{s}$
and the ratio of every pair of unequal diagonal entries of~%the matrix
$c$
is not a~degree~%
$s$
root of unity.
The equality
\[
\left(
\begin{array}{cc}
{c} & {d}
\\
{0} & {\delta }
\end{array}
\right)^{s}
=
\left(
\begin{array}{cc}
{c^{s} } & {(c^{s-1} +\delta c^{s-2} +\dots+\delta ^{s-2} c+\delta
  ^{s-1} e)d}
\\
{0} & {\delta ^{s} }
\end{array}
\right),
\]
where~%
$e$
is the identity matrix of size~%
$n$,
shows that
we can obtain every matrix in
$v(T_{n+1} (K),x^{s} )$
by raising to power~%
$s$
some matrix in
$T_{n+1} (K)$
found by choosing a~suitable vector~%
$d$
provided that
$f=c^{s-1} +\delta c^{s-2} +\dots+\delta ^{s-2} c+\delta ^{s-1} e$
is an~invertible matrix.
However,
if~%the matrix
$f$
is not invertible
then it has a~zero main diagonal entry.
As in the case
$n=2$,
if at least two diagonal entries of~$f$
vanish
then by the inductive assumption
all corresponding diagonal entries of~$c$
are equal.
Indeed,
if
$f_{ii} = f_{jj} = 0$
with
$i\ne j$
then
$c_{ii}/\delta$
and
$c_{jj}/\delta$
are degree~%
$s$
roots of unity,
and so is
$c_{ii}/c_{jj}$.
Therefore,
$c_{ii}=c_{jj}$.
Replacing~%
$\delta$
with
$c_{ii}$,
we obtain an~invertible matrix~%
$f$.

(b)
Suppose now that
the exponent~%
$s$
is divisible by the characteristic~%
$p$
of the field~%
$K$
and,
in the case of finite~%
$K$,
that
%the exponent~%
$s$
is not divisible by
$\left|K\right|-1$.
The relation
\[
\left(
\begin{array}{cc}
{\alpha } & {\beta }
\\
{0} & {\alpha }
\end{array}
\right)^{p}
=
\left(
\begin{array}{cc}
{\alpha ^{p} } & {p\alpha ^{p-1} \beta }
\\
{0} & {\alpha ^{p} }
\end{array}
\right)
=
\left(
\begin{array}{cc}
{\alpha ^{p} } & {0}
\\
{0} & {\alpha ^{p} }
\end{array}
\right),
\]
the uniqueness of degree~%
$p$
roots in~%the field
$K$,
and the description of the set of verbal subgroups of
the groups of triangular matrices
show that
the width of the verbal subgroup
$v(T_{2} (K),x^{p} )$
with respect to the set of generators
equal to the
$p$th
powers of the matrices in
$T_{2} (K)$
is greater than~1.
Similar arguments are valid for every~%
$s$
dividing
$p$
and every
$n\ge 2$.

Take the identity matrix~%
$e$
and an~arbitrary matrix~%
$a$
in
$T_{n-1} (K)$,
column vectors
$b$~and~$c$
of height
$n-1$,
and nonzero elements
$\alpha$~and~$\beta$
of the field~%
$K$.
If a~nonzero element~%
$\beta$
of~%
$K$
is not a~degree~%
$s$
root of unity
then the relation
\[
\left(
\begin{array}{cc}
{\beta ^{-1} a} & {0}
\\
{0} & {\alpha }
\end{array}
\right)^{s}
\left(
\begin{array}{cc}
{\beta e} & {c}
\\
{0} & {1}
\end{array}
\right)^{s}
=
\left(
\begin{array}{cc}
{a^{s} } & {(\beta ^{s-1} +\beta ^{s-2} +\dots+\beta +1)\beta ^{-s} a^s c}
\\
{0} & {\alpha ^{s} }
\end{array}
\right)
\]
shows that
by choosing a~suitable vector~%
$c$
we can express an~arbitrary matrix
$$
\left(
\begin{array}{cc}
{a^{s} } & {b}
\\
{0} & {\alpha ^{s}}
\end{array}
\right)
$$
in
$v(T_{n} (K),x^{s} )$
as the product of two matrices
which are the
$s$th
powers of matrices in
$T_{n} (K)$.
Consequently,
claim~(b) of Theorem~1 is established as well.

(c)
Since in this case
the multiplicative group of the field~%
$K$
is a~cyclic group of order
$\left|K\right|-1$,
it follows that
$a^{\left|K\right|-1}$
lies in
$UT_{n} (K)$
for every
$a \in T_{n} (K)$.
Since
$UT_{n} (K)$
is a~%
$p$-group
and
$(\left|K\right|-1)^{u} r$
is coprime to~%
$p$,
every matrix in
$UT_{n} (K)$
equals itself raised to power
$(\left|K\right|-1)^{u} r$.
Consequently,
$v(T_{n} (K),x^{(\left|K\right|-1)^{u} r} )=UT_{n} (K)$
and
$\wid(T_{n} (K),x^{(\left|K\right|-1)^{u} r} )=1$.
Then
$$
v(T_{n} (K),x^{s} )
=v(v(T_{n} (K),x^{(\left|K\right|-1)^{u} r}),x^{p^{t} } )
=v(UT_{n} (K),x^{p^{t} } ).
$$
It is shown in~\cite{5} that
$\wid(UT_{n} (K),x^{p^{s} } )=1^{}$
for
$n\le p^{t} +2$
and
$\wid(UT_{n} (K),x^{p^{t} } )=2$
for
$n\ge p^{t} +3$.
This suffices to complete the proof of claim~(c) of Theorem~1.
\end{proof}

Since the group
$FT(K)$
consists of the invertible infinite triangular matrices over the field~%
$K$
which differ from the identity matrix
only in finitely many entries,
Theorem~2 easily follows from Theorem~1.

\begin{proof}[Proof of Theorem~3]
Recall that
for every field~%
$K$
with more than two elements
and all positive integers
$r$~and~$s$
we have %the relations
\begin{align*}
\left[ UT_{n}^{r} (K),\; T_{n} (K) \right]
&=
UT_{n}^{r} (K),
\\
\left[UT_{n}^{r} (K),\; UT_{n}^{s} (K)\right]
&=
UT_{n}^{r+s} (K),
\end{align*}
where
$UT_{n}^{r} (K)$
is a~subgroup of %the group
$UT_{n} (K)$
consisting of the matrices with
$r-1$
zero diagonals above the main diagonal
(see [4, p.~40] for instance).
Thus,
in order to prove Theorem~3
by induction on the number of commutators in the word~%
$w$
it suffices to show that
every matrix in the group
appearing on the right-hand side of each of these relations
is the commutator of two matrices in the groups
appearing on the left-hand sides.
Consider two cases corresponding to these relations.

1.
Take an~arbitrary matrix~%
$c$
in
$UT_{n}^{r} (K)$
and find a~matrix~%
$a$
in
$UT_{n}^{r} (K)$
and a~matrix~%
$b$
in
$T_{n} (K)$
so that
$\left[a,b\right]=c$.
Since the field~%
$K$
contains at least three elements,
we can choose a~matrix
$a=e+\sum _{i=1}^{n-r}a_{i,i+r} e_{i,i+r}$
satisfying %the inequalities
$a_{i,i+r} \ne 0$
and
$a_{i,i+r} +c_{i,i+r} \ne 0$
for
$i=1,\dots,n-r$.
Here,
as usual,
$e$
stands for the identity matrix,
and
$e_{i,i+r}$
is a~matrix unit of size~%
$n$.
Expanding the matrix equality
$ab=bac$
entrywise,
we obtain %the equality
\begin{align}
a_{i,i+r} b_{i+r,i+r}
&=
b_{i,i} (a_{i,i+r} +c_{i,i+r} ), \quad
i=1,2,\dots,n-r,
\\
a_{i,i+r} b_{i+r,j+r}
&=
f_{i,j}, \quad
i=1,2,\dots,n-r,\;
j=i+1,\dots,n-r,
\end{align}
in~%
$K$,
where
$f_{i,j}$
are linear functions of the entries in row~%
$i$
of~% the matrix
$b$
with coefficients
obtained from the entries of %the matrices
$a$~and~$c$
using addition and multiplication.
If we assign arbitrary values in~%
$K$
to the entries in rows~1 through~$r$
of~%
$b$
so that
the diagonal entries are nonzero
then,
using~%the equality
(1),
we can determine the diagonal entries
$b_{i+1,i+1} ,\dots,b_{n,n}$
of~%the matrix
$b$,
which moreover are all nonzero,
and using~%the equality
(2),
we can successively determine
the values of the off-diagonal entries in rows
$r+1$
through~%
$n$
of~%the matrix
$b$.

2.
Take now an~arbitrary matrix~%
$c$
in
$UT_{n}^{r+s} (K)$
and find a~matrix~%
$a$
in
$UT_{n}^{r} (K)$
and a~matrix~%
$b$
in
$UT_{n}^{s} (K)$
satisfying
$[a,b]=c$.
%As~%
%$a$
%we take
Put
$a=e+\sum _{i=1}^{n-r}e_{i,i+r}$.
Expanding the equality
$ab=bac$
entrywise,
we obtain %the equality
\begin{equation}
b_{i+r,j+r+s} =f_{i,j} , \quad
i=1,\dots,n-r-s,\;
j=i+1,\dots,n-r-s,
\end{equation}
where
$f_{i,j}$
are linear functions of the entries in row~%
$i$
of~%the matrix
$b$
with coefficients obtained from the entries of~%the matrix
$c$
using addition and multiplication.
If we assign arbitrary values to the off-diagonal entries in rows~%
$1$
through~%
$r$
of~%the matrix
$b$
then %the equalities
(3)
enables us to successively determine the entries in rows
$(r+1)$
through~%
$n$
of~%the matrix
$b$.
\end{proof}

Theorem~4 follows easily from Theorem~3.

\end{document}